\documentclass{article}
\usepackage{amsmath,amssymb}
\newtheorem{df}{Definition}[section]
\newtheorem{thm}[df]{Theorem}
\newtheorem{lem}[df]{Lemma}
\newtheorem{prop}[df]{Proposition}
\newtheorem{cor}[df]{Corollary}
\title{Poisson wavelets on $n$-dimensional spheres}
\author{Ilona Iglewska-Nowak\footnote{West Pomeranian University of Technology, School of Mathematics, al. Piast\'ow 17, 70--310 Szczecin, Poland}}

\begin{document}

\maketitle

\bibliographystyle{amsplain}

\begin{abstract}In this paper, Poisson wavelets on $n$-dimensional spheres, derived from Poisson kernel, are introduced and characterized. We compute their Gegenbauer expansion with respect to the origin of the sphere, as well as with respect to the field source. Further, we give recursive formulae for their explicit representations and we show how the wavelets are localized in space. Also their Euclidean limit is calculated explicitly and its space localization is described. We show that Poisson wavelets can be treated as wavelets derived from approximate identities and we give two inversion formulae.
\end{abstract}

\begin{bfseries}Keywords:\end{bfseries} Poisson wavelets, $n$-spheres \\
\begin{bfseries}AMS Classification:\end{bfseries} 42C40

\section{Introduction}
Investigation of data on higher dimensional spheres has become more and more important in the last years. Recently, a nonzonal continuous wavelet transform for $n$--spheres was introduced~\cite{EBCK09}, being a generalization of wavelet transforms for the two--sphere~\cite{FW96,FW-C,FGS-book} and for the three--sphere~\cite{sB09,BE10KBW,BE10WS3} (both based on the theory of singular integrals~\cite{LV}) to $n$--dimensional nonzonal case. Further development of this theory is to be found in~\cite{IIN13}, and a special case of diffusive wavelets is studied in~\cite{sE11}. Only one different approach for continuous wavelet transform for $n$--spheres is known to the author, namely the wavelet transform based on group--theoretical considerations presented in~\cite{AVn}. 

Motivated by the fact that Poisson wavelets on the two--dimensional sphere~\cite{HI07} have shown to be useful in computations~\cite{HCM03}, and that continuous wavelet transform with respect to them can be discretized and yields a frame~\cite{IH10}, in the present paper, we generalize this notion to $n$ dimensions. In~\cite{IIN13b} we show that the wavelet transform can be discretized to a wavelet frame.

The paper is organized as follows. Section~\ref{sec:sphere} contains basic information about analysis of functions on spheres. In Section~\ref{sec:def} we introduce Poisson wavelets as derivatives of Poisson kernel located inside the unit ball, and we state their basic properties, such as Gegenbauer expansion with respect to the origin of the coordinate system. The Gegenbauer expansion of harmonic continuation of Poisson wavelets with respect to the location of the kernel is derived in Section~\ref{sec:harmonic_continuation}. In Section~\ref{sec:explicit} we give explicit expressions for Poisson wavelets as irrational functions of the first spherical coordinate~$\theta_1$ (since they are zonal functions) and the distance~$r$ of the source location from the origin of coordinate system. Using them, in Section~\ref{sec:localization} we compute the space localization of wavelets, and using the representation from Section~\ref{sec:harmonic_continuation} we find explicit expressions for the Euclidean limit in Section~\ref{sec:euclidean}. In the end, we show that Poisson wavelets are bilinear (Section~\ref{sec:bilinear_wv}) and linear (Section~\ref{sec:linear_wv}) wavelets according to definitions given in~\cite{IIN13}, and consequently that they possess all the advantageous properties of these constructions. In particular, we give two inversion formulae for the wavelet transform with respect to Poisson wavelets.

\section{Preliminaries}\label{sec:sphere}

\subsection{Functions on the sphere}

By $\mathcal{S}^n$ we denote the $n$--dimensional unit sphere in $(n+1)$--dimensional Euclidean space~$\mathbb{R}^{n+1}$ with the rotation--invariant measure~$d\sigma$ normalized such that
$$
\Sigma_n=\int_{\mathcal{S}^n}d\sigma=\frac{2\pi^{\lambda+1}}{\Gamma(\lambda+1)},
$$
where~$\lambda$ and~$n$ are related by
$$
\lambda=\frac{n-1}{2}.
$$
The surface element $d\sigma$ is explicitly given by
$$
d\sigma=\sin^{n-1}\theta_1\,\sin^{n-2}\theta_2\dots\sin\theta_{n-1}d\theta_1\,d\theta_2\dots d\theta_{n-1}d\varphi,
$$
where $(\theta_1,\theta_2,\dots,\theta_{n-1},\varphi)\in[0,\pi]^{n-1}\times[0,2\pi)$ are spherical coordinates satisfying
\begin{align*}
x_1&=\cos\theta_1,\\
x_2&=\sin\theta_1\cos\theta_2,\\
x_3&=\sin\theta_1\sin\theta_2\cos\theta_3,\\
&\dots\\
x_{n-1}&=\sin\theta_1\sin\theta_2\dots\sin\theta_{n-2}\cos\theta_{n-1},\\
x_n&=\sin\theta_1\sin\theta_2\dots\sin\theta_{n-2}\sin\theta_{n-1}\cos\varphi,\\
x_{n+1}&=\sin\theta_1\sin\theta_2\dots\sin\theta_{n-2}\sin\theta_{n-1}\sin\varphi.
\end{align*}
$\left<x,y\right>$ or $x\cdot y$ stand for the scalar product of vectors with origin in~$O$ and endpoints on the sphere. As long as it does not lead to misunderstandings, we identify these vectors with points on the sphere.

Scalar product of $f,g\in\mathcal L^2(\mathcal S^n)$ is defined by
$$
\left<f,g\right>_{\mathcal L^2(\mathcal S^n)}=\frac{1}{\Sigma_n}\int_{\mathcal S^n}\overline{f(x)}\,g(x)\,d\sigma(x),
$$
and by $\|\circ\|$ we denote the induced $\mathcal L^2$--norm.

Gegenbauer polynomials $C_l^\lambda$ of order~$\lambda$, $\lambda\in\mathbb R$, and degree~$l\in\mathbb{N}_0$, are defined in terms of their generating function
\begin{equation}\label{eq:gen_fct_Gg}
\sum_{l=0}^\infty C_l^\lambda(t)\,r^l=\frac{1}{(1-2tr+r^2)^\lambda},\qquad t\in[-1,1],
\end{equation}
and they are explicitly given by
\begin{equation}\label{eq:Gegenbauer_explicitly}
C_l^\lambda(t)=\sum_{k=0}^{[l/2]} (-1)^k\frac{\Gamma(l-k+\lambda)}{\Gamma(\lambda)\,k!\,(l-2k)!}\,(2t)^{l-2k},
\end{equation}
cf. \cite[Sec. IX.3.1, formula (3)]{Vilenkin}.
They are real-valued and for some fixed $\lambda\ne0$ orthogonal to each other with respect to the weight function~$\left(1-\circ^2\right)^{\lambda-\frac{1}{2}}$, compare~\cite[formula~8.939.8]{GR}.

Let $Q_l$ denote a polynomial on~$\mathbb{R}^{n+1}$ homogeneous of degree~$l$, i.e., such that $Q_l(az)=a^lQ_l(z)$ for all $a\in\mathbb R$ and $z\in\mathbb R^{n+1}$, and harmonic in~$\mathbb{R}^{n+1}$, i.e., satisfying $\Delta Q_l(z)=0$, then $Y_l(x)=Q_l(x)$, $x\in\mathcal S^n$, is called a hyperspherical harmonic of degree~$l$. The set of hyperspherical harmonics of degree~$l$ restricted to~$\mathcal S^n$ is denoted by $\mathcal H_l=\mathcal H_l(\mathcal S^n)$. $\mathcal H_l$--functions are eigenfunctions of Laplace--Beltrami operator $\Delta^\ast:=\left.\Delta\right|_{\mathcal S^n}$ with eigenvalue $-l(l+2\lambda)$, further, hyperspherical harmonics of distinct degrees are orthogonal to each other. The number of linearly independent hyperspherical harmonics of degree~$l$ is equal to
$$
N=N(n,l)=\frac{(n+2l-1)(n+l-2)!}{(n-1)!\,l!}.
$$

Every $\mathcal{L}^1(\mathcal S^n)$--function~$f$ can be expanded into Laplace series of hyperspherical harmonics by
$$
S(f;x)\sim\sum_{l=0}^\infty Y_l(f;x),
$$
where $Y_l(f;x)$ is given by
$$
Y_l(f;x)=\frac{\Gamma(\lambda)(\lambda+l)}{2\pi^{\lambda+1}}\int_{\mathcal S^n}C_l^\lambda(x\cdot y)\,f(y)\,d\sigma(y)
   =\frac{\lambda+l}{\lambda}\left<C_l^\lambda(x\cdot\circ),f\right>.
$$
For zonal functions (i.e., those depending only on~$\theta_1=\left<\hat e,x\right>$, where~$\hat e$ is the north-pole of the sphere
$\hat e=(1,0,\dots,0)$) we obtain the Gegenbauer expansion
$$
f(\cos\theta_1)=\sum_{l=0}^\infty\hat f(l)\,C_l^\lambda(\cos\theta_1)
$$
with Gegenbauer coefficients
$$
\hat f(l)=c(l,\lambda)\int_{-1}^1 f(t)\,C_l^\lambda(t)\left(1-t^2\right)^{\lambda-1/2}dt,
$$
where~$c$ is a constant that depends on~$l$ and~$\lambda$. Such functions are identified with functions over the interval $[-1,1]$, i.e., whenever it does not lead to mistakes, we write
$$
f(x)=f(\cos\theta_1).
$$

For $f,g\in\mathcal L^1(\mathcal S^n)$, $g$ zonal, their convolution $f\ast g$ is defined by
\begin{equation*}
(f\ast g)(x)=\frac{1}{\Sigma_n}\int_{\mathcal S^n}f(y)\,g(x\cdot y)\,d\sigma(y).
\end{equation*}
With this notation we have
$$
Y_l(f;x)=\frac{\lambda+l}{\lambda}\,\left(f\ast C_l^\lambda\right)(x),
$$
i.e., the function
$$
\mathcal K_l^\lambda=\frac{\lambda+l}{\lambda}\,C_l^\lambda
$$
is the reproducing kernel for~$\mathcal H_l$.

\section{Definition and basic properties}\label{sec:def}
In the next three sections, we introduce Poisson wavelets on $n$--spheres and prove their basic properties. It is a generalization od ideas from~\cite{HCM03,HI07,IIN07} to the case of $n$--dimensional spheres. Our motivation is the fact that zonal Poisson wavelets on~$\mathcal S^2$ and also their directional counterpart have proven to be of great practical importance and very handful for implementation~\cite{HCM03,CPMDHJ05,HH09}. Since in the recent years investigation of objects in higher dimensions becomes popular (e.g., in crystallography~\cite{sB09,BE10KBW,BE10WS3}), we want to propose an ananalytic tool for their analysis.

Let the Poisson kernel for the unit sphere be given,
\begin{equation}\label{eq:Poisson_kernel}
p_\zeta(y)=\frac{1}{\Sigma_n}\frac{1-|\zeta|^2}{|\zeta-y|^{n+1}}=\frac{1}{\Sigma_n}\frac{1-r^2}{(1-2r\cos\theta+r^2)^{(n+1)/2}},
\end{equation}
where $\zeta,\,y\in\mathbb R^{n+1}$,
$$
r=|\zeta|<|y|=1
$$
and $\theta$ is the angle between the vectors~$\zeta$ and~$y$, i.e.,
$$
r\cos\theta=\zeta\cdot y,
$$
compare~\cite{SW71}. Without loss of generality, suppose $\zeta$ lies on the positive $x_1$--axis, i.e.,
$$
\zeta=r\hat e.
$$
In this case, $\theta$ is the $\theta_1$--coordinate of~$y$. Since
$$
\frac{1-r^2}{(1-2r\cos\theta+r^2)^{\lambda+1}}=\left(1+\frac{r}{\lambda}\,\partial_r\right)\frac{1}{(1-2r\cos\theta+r^2)^\lambda},
$$
we obtain by the generating function equation~\eqref{eq:gen_fct_Gg}
\begin{equation}\label{eq:PKseries}
p_\zeta(y)=\frac{1}{\Sigma_n}\,\sum_{l=0}^\infty r^l\,\mathcal K_l^\lambda(\cos\theta),
\end{equation}
compare also Theorem~IV.2.10 and Theorem~IV.2.14 in~\cite{SW71}. Note that
$$
\Psi_\zeta(x):=\frac{1}{\Sigma_n}\frac{1}{(|x|^2-2x\cdot\zeta+|\zeta|^2)^\lambda},\qquad x\in\mathbb R^{n+1}\setminus\{\zeta\},
$$
is the field caused by a monopole inside the unit ball,
$$
\Delta \Psi_\zeta=\delta_{\zeta},
$$
where~$\delta_\zeta$ is the Dirac measure located at~$\zeta$, consequently,
\begin{equation}\label{eq:multipoles}
p_\zeta=\Psi_\zeta+\frac{1}{\lambda}\,\Psi_\zeta^1,
\end{equation}
where~$\Psi_\zeta^m$ denotes the field caused by a multipole,
$$
\Psi_\zeta^m=(r\partial_r)^m\Psi_\zeta=\frac{1}{\Sigma_n}\sum_{l=0}^\infty l^mr^l C_l^\lambda,\qquad\Delta\Psi_\zeta^m=(r\partial_r)^m\delta_{\zeta},\qquad\zeta=r\hat e.
$$

\begin{df}Poisson wavelet of order~$m$, $m\in\mathbb N$, at a scale~$a$, $a\in\mathbb R_+$, is given recursively by
\begin{align}
g_a^1&=ar\partial_rp_{r\hat e},\qquad r=e^{-a},\notag\\
g_a^{m+1}&=ar\partial_rg_a^m.\label{eq:recursion}
\end{align}
\end{df}

\begin{lem}Gegenbauer expansion of Poisson wavelets is given by
$$
g_a^{m}(y)=\frac{1}{\Sigma_n}\,\sum_{l=0}^\infty\frac{\lambda+l}{\lambda}\,(al)^m e^{-al}\,\mathcal C_l^\lambda(\cos\theta).
$$
\end{lem}
\begin{bfseries}Proof. \end{bfseries} Use formula~\eqref{eq:PKseries}.

\begin{lem}Poisson wavelet of order~$m$ is a sum of fields caused by multipoles,
\begin{equation}\label{eq:sum_of_multipoles}
g_a^m=a^m\left(\Psi_{r\hat e}^m+\frac{1}{\lambda}\,\Psi_{r\hat e}^{m+1}\right)
\end{equation}
\end{lem}
\begin{bfseries}Proof. \end{bfseries}Apply $m$ times the operator $(r\partial_r)$ to the equation~\eqref{eq:multipoles}.

\section{Harmonic continuation}\label{sec:harmonic_continuation}

A remarkable feature of Poisson wavelets is that they possess a representation as a finite sum of hyperspherical harmonics, however, such centered in the point where the field source (multipole) is located. Moreover, a harmonic continuation exists to functions over the space with the source point excluded.

\begin{prop} Poisson wavelets~$g_a^m$, $m\in\mathbb N$, can be uniquely harmonically continued to functions over $\mathbb R^{n+1}\setminus\{r\hat e\}$. They are given by
\begin{equation}\label{eq:finite_multipoles}
g_a^m(x)=\frac{a^m}{\Sigma_n}\sum_{l=0}^{m+1}l!\left(\alpha_l^m+\frac{\alpha_l^{m+1}}{\lambda}\right)e^{-al}\,\frac{C_l^\lambda(\cos\chi)}{|x-r\hat{e}|^{l+2\lambda}},
\end{equation}
where
$$
\cos\chi=\frac{x-r\hat{e}}{|x-r\hat{e}|}\cdot\hat{e}
$$
and the coefficients~$\alpha_l^m$ are recursively given by
\begin{align*}
\alpha_0^0&=1,\\
\alpha_0^m&=0\qquad\text{for }m\geq1,\\
\alpha_m^l&=0\qquad\text{for }l>m,\\
\alpha_l^{m+1}&=l\alpha_l^m+\alpha_{l-1}^m.
\end{align*}
\end{prop}

The proof is analogous to the proof of~\cite[Proposition~1]{HI07}.

\begin{bfseries}Proof. \end{bfseries}Let $x\in\mathbb R^{n+1}$ be different from~$r\hat e$. For its distance from~$r\hat e$ we obtain the expression
\begin{equation*}\begin{split}
|x-r\hat{e}|^2&=|x|^2\left|\hat{x}-\frac{r}{|x|}\,\hat{e}\right|^2
      =|x|^2\left(\sin^2\theta+\left(\cos\theta-\frac{r}{|x|}\right)^2\right)\\
      &=|x|^2\left(1-\frac{2r}{|x|}\,\cos\theta+\left(\frac{r}{|x|}\right)^2\right)
\end{split}\end{equation*}
for $\hat{x}=\frac{x}{|x|}$ and $\cos\theta=\hat{x}\cdot\hat{e}$.
According to~\eqref{eq:gen_fct_Gg} we have
\begin{equation}\label{eq:Psi_finite_sum}
\Sigma_n\Psi_{r\hat e}(x)=\frac{1}{|x|^{2\lambda}}\sum_{l=0}^\infty\left(\frac{r}{|x|}\right)^lC_l^\lambda(\cos\theta),
\end{equation}
and therefore
$$
\left.\Sigma_n\partial_r^m\Psi_{r\hat e}(x)\right|_{r=0}=m!\,\frac{C_m^\lambda(\cos\theta)}{|x|^{m+2\lambda}}.
$$
Consequently, for a field caused by a multipole located at $\zeta=r\hat e$ we may write
$$
\Sigma_n\partial_r^m\Psi_{r\hat e}(x)=m!\,\frac{C_m^\lambda(\cos\chi)}{|x-r\hat e|^{m+2\lambda}}.
$$
We put this expression into~\eqref{eq:sum_of_multipoles} and obtain the representation~\eqref{eq:finite_multipoles} for the wavelet~$g_a^m$, where the coefficients~$\alpha_l^m$ are defined through
\begin{equation}\label{eq:coeffs_alpha}
(r\partial_r)^m=\sum_{l=1}^{m+1}\alpha_l^mr^l\partial_r\text{ for }l\leq m
\end{equation}
and
$$
\alpha_l^m=0\text{ for }l>m.
$$
The recursive formula follows from~\eqref{eq:coeffs_alpha}.\hfill$\Box$

\begin{prop} The harmonically extended Poisson wavelet has the following expansion around~$0$:
$$
g_a^m(y)=\frac{a^m}{\Sigma_n|x|^{2\lambda}}\sum_{l=0}^\infty l^m\left(\frac{r}{|x|}\right)^l\mathcal K_l^\lambda(\cos\theta),
$$
$y\in\mathcal S^n$, $x\in\mathbb R^{n+1}$, $r=e^{-a}$, $\cos\theta=\left<\hat e,y\right>$.
\end{prop}

\begin{bfseries}Proof. \end{bfseries}Apply the operator
$$
(r\partial _r)^m+\frac{(r\partial _r)^{m+1}}{\lambda}
$$
to the equation~\eqref{eq:Psi_finite_sum}.\hfill$\Box$

\section{Explicit expressions}\label{sec:explicit}
In this section we derive explicit formulae for Poisson wavelets as irrational functions of the first spherical variable. This is one of the features that makes them suited for applications. Note that Gauss--Weierstrass wavelets~\cite[Sec.~10]{FGS-book} as well as all the discrete wavelets investigated in~\cite[Sec.~11]{FGS-book} are given only as Laplace series.

\begin{prop} Poisson wavelets of order $m\in\mathbb N$ are represented by
\begin{equation}\label{eq:gam_polynomial}
g_a^m(y)=\frac{a^m}{\Sigma_n}\,D_{\lambda+m+1}\sum_{k=0}^mR_k^{m}(r)\cos^k\theta,
\end{equation}
where
$$
D_j=D_j(r,\theta)=\frac{r}{(1-2r\cos\theta+r^2)^j}
$$
and~$R_k^{m}$ are polynomials of degree~$2m-k+1$, explicitly given by
$$
R_k^{m}(r)=\sum_{j=0}^{[(2m-k+1)/2]}a_j^{m,k}r^{2j+(k-1)_{\text{mod\,}2}},
$$
where the coefficients~$a_j^{m,k}$ satisfy the recursion
\begin{align*}
a_j^{m+1,0}&=b_j^{m+1,0},&j=0,\dots,m+1,\\[0.5em]
a_j^{m+1,k}&=b_j^{m+1,k}+c_j^{m+1,k},&\genfrac{}{}{0pt}{}{k=1,\dots,m,}{j=0,\dots,m+1-\left[\frac{k-1}{2}\right],}\\[0.5em]
a_j^{m+1,m+1}&=c_j^{m+1,m+1},&j=0,\dots,\left[\frac{m+1}{2}\right],
\end{align*}
with
$$\begin{array}{ll}
a_0^{1,0}=-(n+3),&a_1^{1,0}=n-1,\\
a_0^{1,1}=n+1,&a_1^{1,1}=-(n-3)
\end{array}$$
and
\begin{align*}
b_0^{m+1,k}&=2\,a_0^{m,k},\\
b_j^{m+1,k}&=2\,(j+1)\,a_j^{m,k}+2\,(j-\lambda-m-1)\,a_{j-1}^{m,k},\qquad j=1,\dots,m-k/2,\\
b_{m+1-k/2}^{m+1,k}&=-(2\lambda+k)\,a_{m-k/2}^{m,k},\\
c_0^{m+1,k}&=2\,(\lambda+m)\,a_0^{m,k-1},\\
c_j^{m+1,k}&=2\,(\lambda+m)\,a_j^{m,k-1}-2\cdot2j\,a_j^{m,k-1}\\
&=2\,(\lambda+m-2j)\,a_j^{m,k-1}\qquad j=1,\dots,m+1-k/2,
\end{align*}
for an even~$k$ and
\begin{align*}
b_0^{m+1,k}&=a_0^{m,k},\\
b_j^{m+1,k}&=(2j+1)\,a_j^{m,k}+\left(2\,(j-\lambda-m)-3\right)\,a_{j-1}^{m,k},\qquad j=1,\dots,m-[k/2],\\
b_{m+1-[k/2]}^{m+1,k}&=-(2\lambda+2m+1)\,a_{m-[k/2]}^{m,k}+(2m-k+1)\,a_{m-[k/2]}^{m,k}\\
&=-(2\lambda+k)\,a_{m-[k/2]}^{m,k},\\
c_0^{m+1,k}&=0,\\
c_j^{m+1,k}&=2\,(\lambda+m-2j+1)\,a_{j-1}^{m,k-1},\qquad j=1,\dots,m+1-[k/2],
\end{align*}
for an odd~$k$.
\end{prop}

The proof for the two--dimensional case is given in~\cite{IIN07}.

\begin{bfseries}Proof. \end{bfseries}For $m=1$ we compute the wavelet as $ar\partial_rp_{r\hat e}$ with~$p_{r\hat e}$ given by~\eqref{eq:Poisson_kernel} and obtain
\begin{align*}
g_a^1(y)=&\frac{ar}{\Sigma_n\cdot(1-2r\cos\theta+r^2)^{(n+3)/2}}\\
&\cdot\left[\left(-(n+3)r+(n-1)r^3\right)+\left((n+1)-(n-3)r^2\right)\cos\theta\right].
\end{align*}
Suppose, $g_a^m$ has the representation~\eqref{eq:gam_polynomial}. Since
$$
\partial_r D_j=\frac{\left[1-(2j-1)r^2+2(j-1)r\cos\theta\right]}{r}\cdot D_{j+1},
$$
for $m+1$ we obtain by~\eqref{eq:recursion}
\begin{align*}
\frac{\Sigma_n}{a^{m+1}}\,g_a^{m+1}&=D_{\lambda+m+2}\left[1-(2\lambda+2m+1)r^2+2(\lambda+m)r\cos\theta\right]\sum_{k=0}^mR_k^m(r)\cos^k\theta\\
&+D_{\lambda+m+2}\left[1+r^2-2r\cos\theta\right]r\sum_{k=0}^mR_k^{m\,\prime}(r)\cos^k\theta.
\end{align*}
Collecting the powers of~$\cos\theta$ leads to
\begin{align*}
\frac{\Sigma_n}{a^{m+1}}\,g_a^{m+1}
   &=D_{\lambda+m+2}\left[B_0^{m+1}+\sum_{k=1}^m(B_k^{m+1}+C_k^{m+1})\cos^k\theta+C_{m+1}^{m+1}\cos^{m+1}\theta\right],
\end{align*}
where the coefficients
$$
B_k^{m+1}=B_k^{m+1}(r)=\sum_{j=0}^{m+1-\left[\frac{k-1}{2}\right]}b_j^{m+1,k}r^{2j+(k-1)_{\text{mod}\,\,2}},\qquad k=0,\dots,m,
$$
and
$$
C_k^{m+1}=C_k^{m+1}(r)=\sum_{j=0}^{m+1-\left[\frac{k-1}{2}\right]}c_j^{m+1,k}r^{2j+(k-1)_{\text{mod}\,\,2}},\qquad k=1,\dots,\mbox{m+1},
$$
are given by
\begin{align*}
B_k^{m+1}&=\left(1-(2\lambda+2m+1)r^2\right)R_k^m(r)+(1+r^2)\,rR_k^{m\,\prime}(r),\\
C_k^{m+1}&=2\,(\lambda+m)\,rR_{k-1}^m(r)-2r^2R_{k-1}^{m\,\prime}(r).
\end{align*}
The derivative of~$R_k^m$ is equal to
$$
R_k^{m\,\prime}(r)=\sum_{j=0}^{m-k/2}(2j+1)\,a_j^{m,k}\,r^{2j}
$$
for an even~$k$ and
$$
R_k^{m\,\prime}(r)=\sum_{j=1}^{m-k/2}2j\,a_j^{m,k}\,r^{2j-1},
$$
for an odd~$k$.
Therefore, the for the coefficients of the polynomials~$B_k^{m+1}$ and~$C_k^{m+1}$ we obtain the formulae:
\begin{align*}
b_0^{m+1,k}&=a_0^{m,k}+a_0^{m,k}=2\,a_0^{m,k},\\
b_j^{m+1,k}&=a_j^{m,k}-(2\lambda+2m+1)\,a_{j-1}^{m,k}+(2j+1)\,a_j^{m,k}+(2j-1)\,a_{j-1}^{m,k}\\
&=2\,(j+1)\,a_j^{m,k}+2\,(j-\lambda-m-1)\,a_{j-1}^{m,k},\qquad j=1,\dots,m-k/2,\\
b_{m+1-k/2}^{m+1,k}&=-(2\lambda+2m+1)\,a_{m-k/2}^{m,k}+(2m-k+1)\,a_{m-k/2}^{m,k}\\
&=-(2\lambda+k)\,a_{m-k/2}^{m,k},\\
c_0^{m+1,k}&=2\,(\lambda+m)\,a_0^{m,k-1},\\
c_j^{m+1,k}&=2\,(\lambda+m)\,a_j^{m,k-1}-2\cdot2j\,a_j^{m,k-1}\\
&=2\,(\lambda+m-2j)\,a_j^{m,k-1}\qquad j=1,\dots,m+1-k/2,
\end{align*}
for an even~$k$ and
\begin{align*}
b_0^{m+1,k}&=a_0^{m,k},\\
b_j^{m+1,k}&=a_j^{m,k}-(2\lambda+2m+1)\,a_{j-1}^{m,k}+2j\,a_j^{m,k}+2\,(j-1)\,a_{j-1}^{m,k}\\
&=(2j+1)\,a_j^{m,k}+\left(2\,(j-\lambda-m)-3\right)\,a_{j-1}^{m,k},\qquad j=1,\dots,m-[k/2],\\
b_{m+1-[k/2]}^{m+1,k}&=-(2\lambda+2m+1)\,a_{m-[k/2]}^{m,k}+(2m-k+1)\,a_{m-[k/2]}^{m,k}\\
&=-(2\lambda+k)\,a_{m-[k/2]}^{m,k},\\
c_0^{m+1,k}&=0,\\
c_j^{m+1,k}&=2\,(\lambda+m)\,a_{j-1}^{m,k-1}-2\cdot(2j-1)\,a_{j-1}^{m,k-1}\\
&=2\,(\lambda+m-2j+1)\,a_{j-1}^{m,k-1},\qquad j=1,\dots,m+1-[k/2],
\end{align*}
for an odd~$k$.\hfill$\Box$

\section{Localization of the wavelets}\label{sec:localization}
In order to prove that Poisson wavelets as introduced in Section~\ref{sec:def} are indeed wavelets in the sense of definition from~\cite{EBCK09} or~\cite{IIN13}, we need to show that the wavelets are localized in space. More exactly, the following inequality
$$
a^ng_a^m\left(\cos(a\theta)\right)\leq\frac{\mathfrak c\cdot e^{-a}}{\theta^{m+n}},\qquad\theta\in\left(0,\frac{\pi}{a}\right],
$$
holds uniformly in~$a$ for a constant~$\mathfrak c$. The lemmas, the theorem, and the corollary in this section are analogous and they can be proven along the same lines as those in~\cite[Section~5]{IH10}.

\begin{lem}
\label{lem:rop} Let $Q_m$, $m\in\mathbb{N}$, be a sequence of polynomials in two variables satisfying the recursion
\begin{equation}  \label{eq:recQ}
Q_{m+1}(r,t)=A_m(r,t)\cdot Q_m(r,t)+B(r,t)\cdot\frac{\partial}{\partial r}Q_m(r,t)
\end{equation}
with
\begin{equation}\label{eq:Am}
A_m(r,t)=1-(\alpha+1)r^2+\alpha rt\quad\text{for some positive }\alpha
\end{equation}
and
$$
B(r,t)=(1+r^2-2rt)\,r,
$$
and such that
\begin{equation}  \label{eq:assQ1}
Q_1(1,1)=0,\quad\text{and}\quad\left.\frac{\partial}{\partial r}\,Q_1(r,1)\right|_{r=1}\ne0.
\end{equation}
Then the polynomial $Q_m(1,\cdot)$, $m\geq2$, has an $\left[(m+1)/2\right]$--fold root in~$1$.
\end{lem}

\begin{lem}
\label{lem:envf} Let~$\{f_d\}$ be a family of functions over $(0,1)\times[0,\pi]$ given by
\begin{align*}
f_1(r,\theta)&=\frac{r\,Q_1(r,\cos\theta)}{(1+r^2-2r\cos\theta)^{1+\lambda}}, \\
f_{m+1}(r,\theta)&=r\,\frac{\partial}{\partial r}\,f_m(r,\theta),
\end{align*}
where~$Q_1$ is a polynomial satisfying~\eqref{eq:assQ1} and $\lambda$ is a positive integer or half--integer. Then, for any $k\geq2[m/2]+2\lambda$ there exists a constant~$\mathfrak{c}$ such that
\begin{equation}  \label{eq:envf}
|f_m(r,\theta)|\leq\mathfrak{c}\cdot\frac{r}{\theta^k},\quad\theta\in(0,\pi],
\end{equation}
uniformly in~$r$. For $m\geq2$, the number $2[m/2]+2\lambda$ is the smallest possible exponent~$k$. If~$Q_1(1,y)$ has a simple root in~$1$, then $2\lambda$ is
the smallest possible exponent~$k$ on the right--hand--side of~\eqref{eq:envf} for $m=1$.
\end{lem}

\begin{bfseries}Remark 1. \end{bfseries}For $\lambda=\frac{1}{2}$ the recursion~\eqref{eq:recQ} holds for~$A_m$ given by~\eqref{eq:Am} with $\alpha=2m-1$, i.e., the first statement in the proof of~\cite[Lemma~3]{IH10} is not correct. However, it does not affect the multiplicity of the root of~$E_d(1,\cdot)$ and, consequently, the forthcoming considerations.

\begin{bfseries}Remark 2. \end{bfseries}In the proof of \cite[Lemma~5]{IH10}, the function~$F$ for $\lambda=0$ should be defined by $F(0,\theta)=\theta^k\sin\theta E_d(0,\cos\theta)$ in order to be continuous on $[0,1]\times[0,\pi]$.

\begin{thm}
Let
\begin{equation*}
\Psi_{r\hat e}^m=\frac{1}{\Sigma_n}\sum_{l=0}^\infty l^mr^l\,C_l^\lambda,\qquad m\in\mathbb{N}_0,
\end{equation*}
be the field on the sphere caused by the multipole (monopole for $m=0$) $\mu=(r\partial_r)^m\delta_{r\hat{e}}$. For any $k\geq2\left[\frac{m}{2}\right]+2\lambda$ there exists a constant~$\mathfrak{c}$ such that
$$
|\Psi_{r\hat e}^m(\cos\theta)|\leq\mathfrak{c}\cdot\frac{r}{\theta^k},\quad\theta\in(0,\pi],
$$
uniformly in~$r$. $2\left[\frac{m}{2}\right]+2\lambda$ is the smallest possible exponent on the right--hand--side of this inequality.
\end{thm}

\begin{cor}\label{cor:localization} For any $k\geq2\left[\frac{m+1}{2}\right]+2\lambda$ there exists a constant~$\mathfrak{c}$ such that
$$
|g_a^m(\cos\theta)|\leq\mathfrak{c}\cdot\frac{a^m\,e^{-a}}{\theta^k},\quad\theta\in(0,\pi],
$$
uniformly in~$r$. $2\left[\frac{m+1}{2}\right]+2\lambda$ is the smallest possible exponent on the right--hand--side of this inequality.
\end{cor}

\begin{thm}
\label{thm:scaling} Let~$g_a^m$ be a Poisson wavelet of order~$m$. Then there exists a constant~$\mathfrak{c}$ such that
$$
|a^ng_a^m\left(\cos(a\theta)\right)|\leq\frac{\mathfrak{c}\cdot e^{-a}}{\theta^{m+n}},\quad\theta\in\left(0,\frac{\pi}{a}\right],
$$
uniformly in~$a$. $m+n$ is the largest possible exponent in this inequality.
\end{thm}

\begin{cor}\label{cor:proporcja} The functions $(a,\theta)\mapsto a^ng_a^m(\cos\theta)$ are bounded by $\mathfrak c\cdot e^{-a}$ uniformly in~$\theta$, and $n$ is the smallest possible exponent.
\end{cor}

\section{Euclidean limit}\label{sec:euclidean}
It is easy to verify that normalized Poisson wavelets satisfy conditions of~\cite[Theorem~3.4]{IIN13}. Consequently, the have the Euclidean limit property, i.e.
there exists a function $f:\,[0,\infty)\to\mathbb C$, square integrable with respect to~$x^{n-1}\,dx$, and such that
$$
\lim_{a\to0}\,a^ng_a^m\left(\cos\left(S^{-1}(a\xi)\right)\right)=g^m(|\xi|),
$$
holds point-wise for every $\xi\in\mathbb R^n$, where $S^{-1}$ denotes the inverse stereographic projection.

\begin{thm}
The Euclidean limits of Poisson wavelets are given by
\begin{equation}\label{eq:Eucl_limit}
g^m(|\xi|)=\frac{1}{\Sigma_n\lambda}\,(m+1)!\,\frac{C_{m+1}^\lambda(1/\sqrt{1+|\xi|^2})}{\left(1+|\xi|^2\right)^{(m+n)/2}}.
\end{equation}
\end{thm}

\begin{bfseries}Proof. \end{bfseries}Analogous to the proof of~\cite[Theorem~1]{HI07}.

\begin{prop}
The functions~$g^n$ decay at infinity polynomially with degree
$$
m+n+(m+1)_{\text{mod }2}.
$$
\end{prop}

\begin{bfseries}Proof. \end{bfseries}For an even~$m$ the lowest power of the argument in the polynomial~$C_{m+1}^\lambda$ is equal to~$1$, compare~\eqref{eq:Gegenbauer_explicitly}, consequently the behavior of the numerator in~\eqref{eq:Eucl_limit} for $|\xi|\to\infty$ is governed by~$|\xi|^{-1}$. In the case of an even~$m$, the function~$C_{m+1}^\lambda(1/\sqrt{1+|\circ|^2})$ behaves in infinity like a constant. For an alternative proof compare~\cite[Corollary~1, Section~4]{HI07}\hfill$\Box$

\begin{prop}
The functions $g^m$ are of zero mean with respect to the measure
$$
d\nu(|\xi|)=\frac{4\,(4|\xi|)^{2\lambda}}{(4+|\xi|^2)^{2\lambda+1}}\,d|\xi|.
$$
\end{prop}

\begin{bfseries}Proof. \end{bfseries}Analogous to the proof of~\cite[Proposition~4]{IIN07}.

\section{The bilinear wavelet transform}\label{sec:bilinear_wv}
In this section we prove that normalized Poisson wavelets are indeed wavelets, in the sense of~\cite[Section~5.2]{EBCK09} (with weight $\alpha(\rho)=\frac{1}{\rho}$). Since their Gegenbauer coefficients satisfy
$$
\hat g_a^m(l)=\frac{\lambda+l}{\lambda}\,\psi(al)
$$
for a function~$\psi$, we need to verify whether conditions of the following definition are satisfied.

\begin{df}\label{def:sphwavelet}
A family $$\left\{\Psi_a=\sum_{l=0}^\infty\psi(al)\,\mathcal K_l^\lambda\right\}_{a\in\mathbb{R}_+}$$ of $\mathcal{L}^1(\mathcal S^n)$--functions is called spherical wavelet if it satisfies the following admissibility conditions:
\begin{enumerate}
\item for $l\in\mathbb{N}_0$
\begin{equation}\label{eq:admwv1}
\int_0^\infty\left|\psi(t)\right|^2\,\frac{dt}{t}=1,
\end{equation}
\item for $R\in\mathbb{R}_+$
\begin{equation}\label{eq:admwv2}
\int_{-1}^1\left|\int_R^\infty\bigl(\overline{\Psi_a}\ast\Psi_a\bigr)(t)\,\frac{da}{a}\right|\,\left(1-t^2\right)^{\lambda-1/2}dt\leq\mathfrak c
\end{equation}
with $\mathfrak c$ independent of~$R$.
\end{enumerate}
\end{df}

\begin{thm}\label{thm:FWwavelets}Normalized Poisson wavelets
$$
G_a^m=\sum_{l=0}^\infty\psi_m(al)\,\mathcal K_l^\lambda
$$
with
$$
\psi_m(t)=\frac{2^m}{\sqrt{\Gamma(2m)}}\,t^me^{-t}
$$
are wavelets in the sense of Definition~\ref{def:sphwavelet}.
\end{thm}

In the proof a property of Poisson wavelets is imposed, which we would like to state on its own.

\begin{lem}\label{lem:L1norm}Poisson wavelets of order $m\in\mathbb N_0$ satisfy
$$
\int_{-1}^1\left|g_a^m(t)\right|\,\left(1-t^2\right)^{\lambda-1/2}dt\leq\mathfrak c\,e^{-a}
$$
with a constant~$\mathfrak c$ independent of~$a$.
\end{lem}

\begin{bfseries}Proof.\end{bfseries} We change the variables $t=\cos\theta$ and divide the integration interval into two parts: $0\leq\theta\leq a$ and $a\leq\theta\leq1$. For small values of~$\theta$ we use the estimation from Corollary~\ref{cor:proporcja}, and for big values of~$\theta$ the estimation from Theorem~\ref{thm:scaling}. Since $\sin\theta\leq\theta$, we have
\begin{align*}
\int_0^\pi|g(\cos\theta)|\,\sin^{2\lambda}\theta\,d\theta&\leq\int_0^a\frac{\mathfrak c\,e^{-a}}{a^{2\lambda+1}}\,\theta^{2\lambda}\,d\theta
   +\int_a^1\frac{\mathfrak c\,e^{-a}\,a^{m+2\lambda+1}}{a^{2\lambda+1}\,\theta^{m+2\lambda+1}}\,\theta^{2\lambda}\,d\theta\\
&\leq\mathfrak c\,e^{-a}.
\end{align*}
\hfill$\Box$

\begin{bfseries}Proof of Theorem~\ref{thm:FWwavelets}.\end{bfseries} $\psi_m$ is normalized such that~\eqref{eq:admwv1} holds. In order to prove that~\eqref{eq:admwv2} is satisfied, note that by the reproducing property of~$\mathcal K_l^\lambda$
$$
\overline{G_a^m}\ast G_a^m=\sum_{l=0}^\infty\left|\psi_m(al)\right|^2\mathcal K_l^\lambda,
$$
and further
\begin{equation}\label{eq:approximate_identity}
\int_R^\infty\left(\overline{G_a^m}\ast G_a^m\right)\,\frac{da}{a}=\sum_{l=0}^\infty\varphi_m(Rl)\,\mathcal K_l^\lambda
\end{equation}
with
$$
\varphi_m(Rl)=\int_R^\infty\left|\psi_m(al)\right|^2\frac{da}{a}=W_m(Rl)\,e^{-2Rl},
$$
where~$W_m$ is a polynomial (of degree $2m-1$). Consequently, \eqref{eq:approximate_identity} is a weighted sum of Poisson wavelets at the scale~$2R$, and according to the Lemma~\ref{lem:L1norm}, it is uniformly $\mathcal L^1$--bounded with respect to the weight $(1-t^2)^{\lambda-1/2}$, as required.\hfill$\Box$

\begin{bfseries}Remark.\end{bfseries} $$G_a^m=\frac{2^m\,\Sigma_n}{\sqrt{\Gamma(2m)}}\,g_a^m.$$

Consequently, all the definitions and theorems from~\cite{EBCK09} and~\cite{IIN13} concerning bilinear wavelets apply to Poisson wavelets. In particular, the wavelet transform is given by
$$
\mathcal W_mf(a,x)=\left<G_{a,x}^m,f\right>=f\ast\overline{G_a^m}(x),
$$
where $g_{a,x}^m$ denotes the wavelet rotated to~$x$ and $f\in\mathcal X(\mathcal S^n)$, and it can be inverted in $\mathcal X$--sense by
\begin{equation}\label{eq:bilinear_inversion}
f(x)=\frac{1}{\Sigma_n}\int_0^\infty\int_{\mathcal S^n}\mathcal W_mf(a,y)\,G_{a,y}^m(x)\,\frac{d\sigma(y)\,da}{a},
\end{equation}
where $\mathcal X$ denotes the space of $\mathcal C$ or $\mathcal L^p$, $1\leq p<\infty$. Further, the image of the bilinear wavelet transform with respect to Poisson wavelets is a reproducing kernel Hilbert space with scalar product
$$
\left<F,G\right>_{\mathcal L^2(\mathcal S^n\times\mathbb R_+)}=\frac{1}{\Sigma_n}\int_0^\infty\!\!\!\int_{\mathcal S^n}\overline{F(a,x)}\,G(a,x)\,\frac{d\sigma(x)\,da}{a}
$$
and with reproducing kernel given by
$$
\Pi^m(a,x;b,y)=\frac{\sqrt{\Gamma(4m)}}{\Gamma(2m)}\frac{(ab)^m}{(a+b)^{2m}}\,G_{a+b}^{2m}(x\cdot y).
$$

Note that Abel--Poisson wavelet introduced in~\cite{FW96} can be treated as Poisson wavelet of order~$1/2$.

\section{Linear wavelet transform}\label{sec:linear_wv}

In a similar way we prove that Poisson wavelets yield linear wavelet analysis. Let us recall Definition~4.1 from~\cite{IIN13}.

\begin{df}\label{def:lsw} Let $\{\Psi_a^L\}_{a\in\mathbb R_+}$ be a family of zonal~$\mathcal L^2$--functions such that the following admissibility conditions are satisfied:
\begin{enumerate}
\item for $l\in\mathbb{N}_0$
$$
\int_0^\infty\widehat{\Psi_a^L}(l)\,\frac{da}{a}=\frac{\lambda+l}{\lambda},
$$
\item for $R\in\mathbb{R}_+$
$$
\int_{-1}^1\left|\int_R^\infty\Psi_a^L(t)\,\frac{da}{a}\right|\,\left(1-t^2\right)^{\lambda-1/2}dt\leq\mathfrak c
$$
with $\mathfrak c$ independent of~$R$.
\end{enumerate}
Then $\{\Psi_a^L\}_{a\in\mathbb R_+}$ is called a spherical linear wavelet. The associated wavelet transform is defined by
$$
\mathcal{W}_\Psi^L f(a,x)=\left<\Psi_{a,x}^L,f\right>_{\mathcal{L}^2(\mathcal{S}^n)}=\bigl(f\ast\overline{\Psi_a^L}\bigr)(x).
$$
for all $f\in\mathcal L^2(\mathcal S^n)$.
\end{df}

\begin{thm}Normalized Poisson wavelets
$$
\widetilde G_a^m=\sum_{l=0}^\infty\gamma_m(al)\,\mathcal K_l^\lambda
$$
with
$$
\gamma_m(t)=\frac{1}{\Gamma(m)}\,t^me^{-t}
$$
are linear wavelets in the sense of Definition~\ref{def:lsw}.
\end{thm}
\begin{bfseries}Proof. \end{bfseries}Analogous to the proof of Theorem~\ref{thm:FWwavelets}.

\begin{bfseries}Remark.\end{bfseries} $$\tilde G_a^m=\frac{\Sigma_n}{\sqrt{\Gamma(m)}}\,g_a^m.$$

Note that the wavelet transform is the same as in the bilinear case. The inversion formula, however, is given by
\begin{equation}\label{eq:linear_inversion}
f(x)=\frac{1}{\Sigma_n}\,\int_0^\infty\mathcal W_m^Lf(a,x)\,\frac{da}{a}
\end{equation}
in $\mathcal C$ and $\mathcal L^p$, $1\leq p<\infty$, sense. This means that the wavelet transform with respect to Poisson wavelets defined in Section~\ref{sec:def} can be inverted with any of the formulae~\eqref{eq:bilinear_inversion} respectively~\eqref{eq:linear_inversion} with constant $\frac{1}{\Sigma_n}$ replaced by
$$
\frac{4^m\,\Sigma_n}{\Gamma(2m)}\qquad\text{respectively}\qquad\frac{1}{\Gamma(m)}.
$$

Further, note that Abel--Poisson L--wavelet introduced in~\cite{FW96} is a Poisson wavelet of order~$1$.

\end{document}